\documentclass[12pt,a4paper]{amsart}
\usepackage[italian, UKenglish]{babel}
\usepackage[latin9]{inputenx} 
\usepackage{fontenc}
\usepackage[mathscr]{eucal}
\usepackage{indentfirst}
\usepackage{verbatim}
\usepackage[margin=2.9cm]{geometry}
\usepackage[colorinlistoftodos]{todonotes}
\usepackage{xcolor}
\usepackage[framemethod=tikz]{mdframed} 
\usepackage{amsfonts}
\usepackage{amsthm}
\numberwithin{equation}{section}
\usepackage{amsmath}
\usepackage{amscd}
\usepackage{amssymb}
\usepackage{mathtools}
\usepackage{mathabx}
\usepackage{mathptmx}
\usepackage{hyperref} 
\makeatletter
\def\@tocline#1#2#3#4#5#6#7{\relax
  \ifnum #1>\c@tocdepth 
  \else
    \par \addpenalty\@secpenalty\addvspace{#2}%
    \begingroup \hyphenpenalty\@M
    \@ifempty{#4}{%
      \@tempdima\csname r@tocindent\number#1\endcsname\relax
    }{%
      \@tempdima#4\relax
    }%
    \parindent\z@ \leftskip#3\relax \advance\leftskip\@tempdima\relax
    \rightskip\@pnumwidth plus4em \parfillskip-\@pnumwidth
    #5\leavevmode\hskip-\@tempdima
      \ifcase #1
       \or\or \hskip 1em \or \hskip 2em \else \hskip 3em \fi%
      #6\nobreak\relax
    \dotfill\hbox to\@pnumwidth{\@tocpagenum{#7}}\par
    \nobreak
    \endgroup
  \fi}
\makeatother
\newcommand{\mycomment}[1]{%
}

\newcommand*{\vertiii}[1]{{\left\vert\kern-0.25ex\left\vert\kern-0.25ex\left\vert #1 \right\vert\kern-0.25ex\right\vert\kern-0.25ex\right\vert}}

\newcommand*{\mres}{\mathbin{\vrule height 1.6ex depth 0pt width 0.13ex\vrule height 0.13ex depth 0pt width 1.3ex}}
\makeatletter
\newcommand\mathcircled[1]{%
  \mathpalette\@mathcircled{#1}%
}
\makeatother
\newtheorem{innercustomgeneric}{\customgenericname}
\providecommand{\customgenericname}{}
\newcommand{\newcustomtheorem}[2]{%
  \newenvironment{#1}[1]
  {%
   \renewcommand\customgenericname{#2}%
   \renewcommand\theinnercustomgeneric{##1}%
   \innercustomgeneric
  }
  {\endinnercustomgeneric}
}
\newcustomtheorem{customprop}{Proposition}
\newcustomtheorem{customcor}{Corollary}
\newcustomtheorem{customdefin}{Definition}
 
\allowdisplaybreaks 
\DeclareMathAlphabet{\mathcal}{OMS}{cmsy}{m}{n}
\DeclareMathOperator{\spt}{spt}

\theoremstyle{plain}
\newtheorem{theor}{Theorem}[section] 
\newtheorem{lem}[theor]{Lemma} 

\newtheorem{prop}[theor]{Proposition}
\newtheorem*{cor*}{Corollary}
\newtheorem*{prop*}{Proposition}
\newtheorem*{lem*}{Lemma}
\theoremstyle{definition}

\newtheorem*{defin*}{Definition}

\newtheorem{no}[theor]{Notation}

\begin{document}
\title{A note on the Slicing of $(k+1)$-Currents in the Heisenberg Group $\mathbb{H}^n$ in the case $k=n$}
\author[C. Ackermann]{Colleen Ackermann}
\address{\textsc{Colleen Ackermann, Ph.D.} \newline \hspace*{1em} Montgomery College, Mathematics, Statistics and Data Science Department, USA}
\email{colleen.ackermann@montgomerycollege.edu}
\author[G. Canarecci]{Giovanni Canarecci}
\address{\textsc{Giovanni Canarecci, Ph.D.} \newline \hspace*{1em} Software Engineer, Former University of Helsinki, Finland} 
\email{giovanni.canarecci@alumni.helsinki.fi}
\begin{abstract} 
This paper aims to expand on the open case $k=n$ regarding Proposition 3.6\cite{GCslicing} and hopefully foster curiosity for its resolution.
\end{abstract}
\maketitle
\tableofcontents


\section{Preliminaries}

This note is as a follow-up to the 2023 PhD thesis of Giovanni Canarecci \cite{GCphd}. It is joint work between Giovanni Canarecci and Colleen Ackermann, from a research visit in summer 2022. The trip of Giovanni Canarecci was partially sponsored by the Doctoral School in Natural Sciences of the University of Helsinki.\\\\
Given the brevity of the note, the preliminaries are not meant to be exhaustive and many of the definitions will be assumed. Precise definitions and preliminaries are in article \cite{GCslicing} and properly referenced to the original authors there.\\

\noindent
Proposition 3.6\cite{GCslicing} is the main result of article \cite{GCslicing} and contains the limiting condition $k\neq n$, whose consequences are discussed in the original work. The limitation $k\neq n$ ultimately comes from Lemma 3.11\cite{GCslicing}, where the cases differ substantially. This note is meant to add some steps in the direction for the case $k= n$ and hopefully foster curiosity for its resolution. We start by stating the proposition:

\begin{prop}[3.6\cite{GCslicing}]\label{next3properties}
Consider an open set $U \subseteq \mathbb{H}^n$, $G\in N_{\mathbb{H},k+1}(U)$, $f \in Lip(U,\mathbb{R})$, $t \in \mathbb{R}$ and $k\neq n$. Then the following properties hold:
\begin{enumerate}
  \setcounter{enumi}{3}
\item
$M \left ( \langle G,f,t+ \rangle \right ) \leq    \ Lip(f) \liminf\limits_{h \rightarrow 0+} \frac{1}{h} \mu_G \left (   U \cap \{ t < f < t+h \}  \right )$.
\item
$\int_a^b M \left ( \langle G,f,t+ \rangle \right ) dt \leq    \ Lip(f) \mu_G \left (   U \cap \{a < f < b \}  \right ), \quad a,b \in \mathbb{R}$.
\item
$ \langle  G,f,t+ \rangle \in N_{\mathbb{H},k}(U)  \quad  \text{for a.e. } t$.
\end{enumerate}
\end{prop}

\noindent
Below we give some of the notations and lemmas needed for the proof:

\begin{no}[1.15\cite{GCslicing}]\label{notationL}
We let $L$ denote the operator
\begin{align*}
L:  {\prescript{}{}\bigwedge}^{n-1} \mathfrak{h}_1  \to {\prescript{}{}\bigwedge}^{n+1} \mathfrak{h}_1, \  \beta  \mapsto d \theta \wedge \beta .
\end{align*}
Furthermore we remind the reader that, if $\gamma \in \Omega^{k-1}$, we can consider the equivalence class
$$
  {\prescript{}{}\bigwedge}^k \mathfrak{h}_1  =       \left \{  \beta \in \Omega^k ; \ \beta =0 \ \text{or} \ \beta \wedge \theta \neq 0 \right  \}  \cong  \frac{\Omega^k}{ \{ \gamma \wedge \theta \} }    ,
$$
where we write $\{ \gamma \wedge \theta \} = \{ \gamma \wedge \theta ; \ \gamma \in \Omega^{k-1} \}$ for short. The equivalence is given by $\beta \mapsto ( \beta)_{\vert_{ {\prescript{}{}\bigwedge}^{k} \mathfrak{h}_1 }}$.\\
In particular, $L$ is an isomorphism (see section $2$ in \cite{RUMIN}) and we can denote
$$
\mathcal{L} (\alpha) :=  L^{-1}   \left ( - \left ( d  \alpha \right )_{\vert_{ {\prescript{}{}\bigwedge}^{n+1} \mathfrak{h}_1 }}  \right ) .
$$
%
%
\end{no}

\begin{no}[1.16\cite{GCslicing}]\label{notation:IJ}
We denote by $ [\alpha]_{I^k}$ an element of the quotient $\frac{\Omega^k}{I^k}$ and  $\omega_{\vert_{ J^k }}$ an element of $J^k$ whenever  $\omega \in   \mathcal{D}^{k} (U)$.
\end{no}

\begin{lem}[3.7\cite{GCslicing}]\label{lemmagamma}
Consider an open set $U \subseteq \mathbb{H}^n$, $f \in Lip(U,\mathbb{R})$, $t \in \mathbb{R}$, $h>0$ fixed and $s \in \mathbb{R}$. Then let the function
\begin{equation*}
\gamma_h (s) := \frac{|s-t| - |s-(t+h)| +h}{2h}.
\end{equation*}
One can observe that
$$
\gamma_h \circ f (p) =
\begin{cases}
0 , \quad & f(p)\leq t,\\
\frac{f(p)-t}{h} , \quad & t < f(p) < t+h,\\
1 , \quad & f(p) \geq t+h,
\end{cases}
$$
$$
\gamma_h \circ f \in \ Lip(U,\mathbb{R}) \quad \text{and} \quad Lip(\gamma_h \circ f) \leq \frac{Lip(f)}{h}.
$$
\end{lem}

\begin{lem}[3.9\cite{GCslicing}]\label{lemmag_i--faseII}
Consider an open set $U \subseteq \mathbb{H}^n$, $f \in Lip(U,\mathbb{R})$, $t \in \mathbb{R}$, $h>0$ fixed and the function $\gamma_h$ defined in Lemma \ref{lemmagamma}. Then we can approximate $\gamma_h \circ f$ uniformly by functions $g_i\in C^\infty (U,\mathbb{R})$ (notationally $ g_i \rightrightarrows   \gamma_h \circ f  $),  so that
$$
 \spt d g_i \subseteq  \{ t < f < t+h \}  \quad \text{and} \quad  \lim\limits_{i \to \infty}  Lip(g_i) = Lip(\gamma_h \circ f).
$$
\end{lem}

\begin{lem}[3.11\cite{GCslicing}]\label{L2}
Consider an open set $U \subseteq \mathbb{H}^n$, $G\in \mathcal{D}_{\mathbb{H},k+1}(U)$, $\omega \in \mathcal{D}^k_{\mathbb{H}}(U)$  and the functions $g_i \in C^\infty (U,\mathbb{R})$  defined in Lemma \ref{lemmag_i--faseII}. Also recall notations \ref{notationL} and \ref{notation:IJ}. Then
$$
\left [ \left ( \partial G \right ) \mres{g_i} -  \partial \left ( G \mres{g_i} \right ) \right ]  (\omega)  =
$$
$$
=
\begin{cases}
G   \left ( \left [ d^{(1)} g_i \wedge  \omega \right ]_{I^{k+1}}  \right ) , 
\quad [\omega ]_{I^{k}} \in \mathcal{D}^k_{\mathbb{H}}(U) = \frac{\Omega^k}{I^k},
\quad \text{if }k<n,\\
G  \Big ( 
 d^{(1)} g_i  \wedge \left ( \omega + \mathcal{L}(\omega) \wedge \theta \right )  + 
 d^{(n+1)} \left (
 \left (
 \mathcal{L}( g_i \omega) -   g_i \mathcal{L}( \omega) 
 \right )
\wedge \theta
 \right ) 
\Big ) , \\
\quad \quad  \quad  \quad \quad  \quad  \quad  \quad \quad  \quad  \quad \quad  \quad  \quad \quad  \quad  \quad [\omega ]_{I^{n}} \in \mathcal{D}^n_{\mathbb{H}}(U) = \frac{\Omega^n}{I^n},
\quad \text{if }k=n,\\
G   \left ( \left ( d^{(1)} g_i \wedge  \omega \right )_{\vert_{J^{k+1} }}  \right ) , 
\quad \omega \in \mathcal{D}^k_{\mathbb{H}}(U) = J^k,
\quad \text{if }k>n.
\end{cases}
$$
\end{lem}

\section{The case $k=n$}

\noindent
If $k=n$ we can consider $U \subseteq \mathbb{H}^n$ an open set, $G\in \mathcal{D}_{\mathbb{H},n+1}(U)$,  
 $[\omega ]_{I^{n}} \in \mathcal{D}^n_{\mathbb{H}}(U) = \frac{\Omega^n}{I^n}$,  and the functions $g_i \in C^\infty (U,\mathbb{R})$  defined in Lemma \ref{lemmag_i--faseII}. Then, by Lemma \ref{L2},
\begin{align}\label{newbeginning}
\begin{aligned}
&\left [ \left ( \partial G \right ) \mres{g_i} -  \partial \left ( G \mres{g_i} \right ) \right ] \left  (  [\omega ]_{I^{n}}  \right )  \\
=&
G  \Big (    d^{(1)} g_i  \wedge \left ( \omega + \mathcal{L}(\omega) \wedge \theta \right )  +    d^{(n+1)} \left (   \left (   \mathcal{L}( g_i \omega) -   g_i \mathcal{L}( \omega)    \right )  \wedge \theta   \right )  \Big ) .
\end{aligned}
\end{align}

\noindent
The right hand side can be partially rewritten using the following lemma.

\begin{lem}[3.14\cite{GCslicing}]\label{L3}
Let  $U \subseteq \mathbb{H}^n$ be open,  $\omega \in \Omega^n$   and the functions $g_i \in C^\infty (U,\mathbb{R})$ be as defined in Lemma \ref{lemmag_i--faseII}. Also recall Notation \ref{notationL}.  Then
$$
\mathcal{L}( g_i \omega) -   g_i \mathcal{L}( \omega) 
= L^{-1}   \left ( - \left ( d^{(1)} g_i \wedge  \omega \right )_{\vert_{ {\prescript{}{}\bigwedge}^{n+1} \mathfrak{h}_1 }}  \right ) .
$$
\end{lem}

\noindent
Furthermore, one can observe that the right hand side of equation (\ref{newbeginning}) is not null because of the following lemma. This observation was not needed in the cases of $k\neq n$ as the definition of $\mathcal{D}_{\mathbb{H}^*}(U)$ did not change between $k$ and $k+1$, making the step immediate.

\begin{lem}[3.15\cite{GCslicing}]\label{appartenenza}
Consider an open set $U \subseteq \mathbb{H}^n$, $G\in \mathcal{D}_{\mathbb{H},n+1}(U)$, $\omega \in \Omega^n$  
 and the functions $g_i \in C^\infty (U,\mathbb{R})$  defined in Lemma \ref{lemmag_i--faseII}. Also recall Notation \ref{notationL}. Then
\begin{align}\label{condJ}
    d^{(1)} g_i  \wedge \left ( \omega + \mathcal{L}(\omega) \wedge \theta \right )  +    d^{(n+1)} \left (   \left (   \mathcal{L}( g_i \omega) -   g_i \mathcal{L}( \omega)    \right )  \wedge \theta   \right )  \in J^{n+1}.
\end{align}
\end{lem}

\noindent
The difficulty of the case with $k=n$ lies in obtaining 
 an inequality similar to property (4) in Proposition \ref{next3properties}. This means that, for $G\in N_{\mathbb{H},n+1}(U)$, $f \in Lip(U,\mathbb{R})$ and $t \in \mathbb{R}$, we wish to bound $M \left ( \langle G,f,t+ \rangle \right )$ from above with a quantity including $ \liminf\limits_{h \rightarrow 0+} \frac{1}{h} \mu_G \left (   U \cap \{ t < f < t+h \}  \right )$.
By Lemma 
3.10\cite{GCslicing}, we know
\begin{align*}
M (\langle G,f,t+ \rangle) \leq &  \liminf\limits_{h \rightarrow 0+}  \lim_{i \to \infty}
M  \left ( 
 \left ( \partial G \right ) \mres{g_i} -  \partial \left ( G \mres{g_i} \right ) 
\right ).
\end{align*}
Then, by lemmas \ref{L2} and \ref{L3}, for $\omega \in \Omega^n$,
\begin{align}
\begin{aligned}
\left [ \left ( \partial G \right ) \mres{g_i} -  \partial \left ( G \mres{g_i} \right ) \right ]  ([\omega ]_{I^{n}})  
=&
G  \Bigg (  d^{(1)} g_i  \wedge \left ( \omega + \mathcal{L}(\omega) \wedge \theta \right ) \\
&+   d^{(n+1)} \left (  
L^{-1} \left ( - \left ( d^{(1)} g_i \wedge  \omega \right )_{\vert_{ {\prescript{}{}\bigwedge}^{n+1} \mathfrak{h}_1 }}  \right )
  \wedge \theta \right ) \Bigg ) .
\end{aligned}
\end{align}


\section{New steps on the case $k=n$}

At the end of the latest published development, we were left with:
\begin{align}
\begin{aligned}\label{2piecesT}
\left [ \left ( \partial G \right ) \mres{g_i} -  \partial \left ( G \mres{g_i} \right ) \right ]  ([\omega ]_{I^{n}})  
=&
G  \Bigg (  d^{(1)} g_i  \wedge \left ( \omega + \mathcal{L}(\omega) \wedge \theta \right ) \\
&+   d^{(n+1)} \left (  
L^{-1} \left ( - \left ( d^{(1)} g_i \wedge  \omega \right )_{\vert_{ {\prescript{}{}\bigwedge}^{n+1} \mathfrak{h}_1 }}  \right )
  \wedge \theta \right ) \Bigg ) .
\end{aligned}
\end{align}

\noindent
Inside the parentheses on the right hand side of equation (\ref{2piecesT}), the first component is
\begin{align}\label{firsttermproof}
 d g_i \wedge   \left ( \omega + \mathcal{L}(\omega) \wedge \theta \right )  = d g_i \wedge    \omega  +  d g_i \wedge  \mathcal{L}(\omega) \wedge \theta.
\end{align}

\noindent
Following steps from the proof of Lemma 3.15\cite{GCslicing}, we will rewrite the second term inside the parentheses on the right hand side of equation (\ref{2piecesT}). First consider
\begin{align*}
 d g_i \wedge \omega
=& \sum_{j=1}^{ 2n+1 }  \sum_{1\leq l_1 \leq \dots \leq l_n \leq 2n+1 } 
 W_j g_i  \omega_{ l_1  \dots  l_n} 
 dw_j  \wedge  dw_{l_1} \wedge \dots  \wedge dw_{l_n} 
\end{align*}
where the notation is defined in Notation 1.4\cite{GCslicing}, and
\begin{align*}
 -\left ( d g_i \wedge \omega \right )_{\vert_{{\prescript{}{}\bigwedge}^{n+1} \mathfrak{h}_1 }}
=& - \sum_{j=1}^{ 2n }  \sum_{1\leq l_1 \leq \dots \leq l_n \leq 2n } 
 W_j g_i  \omega_{ l_1  \dots  l_n} 
 dw_j  \wedge  dw_{l_1} \wedge \dots  \wedge dw_{l_n} .
\end{align*}
Notice that the $dw_{l_m}$'s are $n$ different basis elements of $\Omega^1$ and they always have their counterpart $dw_{l_m	\pm n}$ among the $ dw_j $'s, since $j=1,\dots,2n$. Hence we can write
\begin{align*}
 -\left ( d g_i \wedge \omega \right )_{\vert_{{\prescript{}{}\bigwedge}^{n+1} \mathfrak{h}_1 }}
= - \sum_{j=1}^{ n }  dw_j  \wedge  dw_{j+n} \wedge \gamma
= d\theta \wedge \gamma
=\gamma \wedge d\theta,
\end{align*}
where $\gamma \in  {\prescript{}{}\bigwedge}^{n-1} \mathfrak{h}_1$. It follows that
\begin{align*}
L^{-1} \left ( -\left ( d g_i \wedge \omega \right )_{\vert_{{\prescript{}{}\bigwedge}^{n+1} \mathfrak{h}_1 }} \right )  =\gamma
\end{align*}
and
\begin{align}
\begin{aligned}\label{secondtermproof}
 d \left (  L^{-1} \left ( - \left ( d g_i \wedge  \omega \right )_{\vert_{ {\prescript{}{}\bigwedge}^{n+1} \mathfrak{h}_1 }}  \right )   \wedge \theta \right ) 
&= d \left (  \gamma   \wedge \theta \right ) \\
&= d   \gamma   \wedge \theta   + \gamma   \wedge d \theta \\
&= d   \gamma   \wedge \theta   - \left ( d g_i \wedge  \omega \right )_{\vert_{ {\prescript{}{}\bigwedge}^{n+1} \mathfrak{h}_1 }}.
\end{aligned}
\end{align}

\noindent
Then, by combining equations (\ref{firsttermproof}) and (\ref{secondtermproof}), we get
\begin{align}
&d g_i  \wedge \left ( \omega + \mathcal{L}(\omega) \wedge \theta \right ) 
+   d \left (  L^{-1} \left ( - \left ( d g_i \wedge  \omega \right )_{\vert_{ {\prescript{}{}\bigwedge}^{n+1} \mathfrak{h}_1 }}  \right )    \wedge \theta \right ) \notag\\
&=
 d g_i \wedge    \omega   
+  d g_i \wedge  \mathcal{L}(\omega) \wedge \theta  
+ d   \gamma   \wedge \theta   
- \left ( d g_i \wedge  \omega \right )_{\vert_{ {\prescript{}{}\bigwedge}^{n+1} \mathfrak{h}_1 }}\notag\\
&=
\left ( d g_i \wedge  \omega \right )_{\vert_{ ({\prescript{}{}\bigwedge}^{n+1} \mathfrak{h}_1)^\perp }}
+  d g_i \wedge  \mathcal{L}(\omega) \wedge \theta  
 + d   \gamma   \wedge \theta  \label{firstsecondcombinedequation}\\
&=
\left ( d g_i \wedge  \omega \right )_{\vert_{ ({\prescript{}{}\bigwedge}^{n+1} \mathfrak{h}_1)^\perp }}
+  d g_i \wedge  \mathcal{L}(\omega) \wedge \theta  
 + d   \left (  \mathcal{L}( g_i \omega) -   g_i \mathcal{L}( \omega)  \right )   \wedge \theta \notag\\
&=
\left ( d g_i \wedge  \omega \right )_{\vert_{ ({\prescript{}{}\bigwedge}^{n+1} \mathfrak{h}_1)^\perp }}
 +  \left (  d     \mathcal{L}( g_i \omega)  - g_i d     \mathcal{L}( \omega)  \right )  \wedge \theta. \tag{3.4.1}\label{firstsecondcombined}
\end{align}

\noindent
Notice that, if
$$
\omega = \omega' + \beta \wedge \theta, \quad \text{with } \omega \in {\prescript{}{}\bigwedge}^{n} \mathfrak{h}_1 \text{ and } \beta \in {\prescript{}{}\bigwedge}^{n-1} \mathfrak{h}_1,
$$
then
\begin{align*}
 d g_i \wedge  \omega 
&=    d g_i \wedge \omega' +  d g_i \wedge \beta \wedge \theta\\
&=    \sum_{j=1}^{2n+1} W_{j} g_i  \theta_j  \wedge \omega' +  d g_i \wedge \beta \wedge \theta.
\end{align*}
Hence, $\left ( d g_i \wedge  \omega \right )_{\vert_{ ({\prescript{}{}\bigwedge}^{n+1} \mathfrak{h}_1)^\perp }}$ from line (\ref{firstsecondcombined}) in equation (\ref{firstsecondcombinedequation}) can be rewritten as
\begin{align*}
\left ( d g_i \wedge  \omega \right )_{\vert_{ ({\prescript{}{}\bigwedge}^{n+1} \mathfrak{h}_1)^\perp }}
&=  W_{2n+1} g_i  \theta  \wedge \omega' +  d g_i \wedge \beta \wedge \theta.
\end{align*}

\noindent
Next we turn our attention to the second term at the end of line (\ref{firstsecondcombined}) in equation (\ref{firstsecondcombinedequation}). We start by showing that
\begin{align*}
d \omega &= d\omega' +d \beta \wedge \theta  +   \beta \wedge  d \theta;\\
- \left (  d \omega \right )_{\vert_{ {\prescript{}{}\bigwedge}^{n+1} \mathfrak{h}_1 }} &=- \left ( d\omega'  \right )_{\vert_{ {\prescript{}{}\bigwedge}^{n+1} \mathfrak{h}_1 }} - \beta \wedge  d \theta;\\
 \mathcal{L}( \omega) &= L^{-1}  \left (  - \left (  d \omega \right )_{\vert_{ {\prescript{}{}\bigwedge}^{n+1} \mathfrak{h}_1 }}   \right )  -  \beta;\\
- g_i d \mathcal{L}( \omega) &=- g_i d L^{-1}  \left (  - \left (  d \omega \right )_{\vert_{ {\prescript{}{}\bigwedge}^{n+1} \mathfrak{h}_1 }}   \right )  + g_i  d\beta.
\end{align*}

\noindent
Furthermore, we have
\begin{align*}
g_i \omega &= g_i \omega' + g_i \beta \wedge \theta;\\
d ( g_i \omega ) &= dg_i \wedge \omega' +g_i d \omega' +d g_i \wedge \beta \wedge \theta  + g_i d  \beta \wedge  \theta   + g_i \beta \wedge  d \theta;\\
- \left (  d g_i \omega \right )_{\vert_{ {\prescript{}{}\bigwedge}^{n+1} \mathfrak{h}_1 }} &= - \left (  dg_i \wedge \omega'  + g_i d \omega' \right )_{\vert_{ {\prescript{}{}\bigwedge}^{n+1} \mathfrak{h}_1 }} -  g_i \beta \wedge  d \theta;\\
 \mathcal{L}( g_i \omega) &= L^{-1}  \left (   - \left (  dg_i \wedge \omega'  + g_i d \omega' \right )_{\vert_{ {\prescript{}{}\bigwedge}^{n+1} \mathfrak{h}_1 }}   \right )  -  g_i \beta;\\
d \mathcal{L}( g_i \omega) &=d  L^{-1}  \left (   - \left (  dg_i \wedge \omega'  + g_i d \omega' \right )_{\vert_{ {\prescript{}{}\bigwedge}^{n+1} \mathfrak{h}_1 }}   \right )  - d g_i \wedge \beta  -  g_i  d \beta .
\end{align*}

\noindent
Then we see that
\begin{align*}
&d \mathcal{L}( g_i \omega) - g_i d \mathcal{L}( \omega) \\
&= d  L^{-1}  \left (   - \left (  dg_i \wedge \omega'  + g_i d \omega' \right )_{\vert_{ {\prescript{}{}\bigwedge}^{n+1} \mathfrak{h}_1 }}   \right )   - g_i d L^{-1}  \left (  - \left (  d \omega \right )_{\vert_{ {\prescript{}{}\bigwedge}^{n+1} \mathfrak{h}_1 }}   \right )  - d g_i \wedge \beta  .
\end{align*}

\noindent
Returning to line (\ref{firstsecondcombined}) in equation (\ref{firstsecondcombinedequation}), we now have
\begin{align*}
&\left ( d g_i \wedge  \omega \right )_{\vert_{ ({\prescript{}{}\bigwedge}^{n+1} \mathfrak{h}_1)^\perp }}
 +  \left (  d     \mathcal{L}( g_i \omega)  - g_i d     \mathcal{L}( \omega)  \right )  \wedge \theta\\
&=
 \bigg [
- W_{2n+1} g_i   \omega' +  d g_i \wedge \beta 
+
d  L^{-1}  \left (   - \left (  dg_i \wedge \omega'  + g_i d \omega' \right )_{\vert_{ {\prescript{}{}\bigwedge}^{n+1} \mathfrak{h}_1 }}   \right ) \\
&\quad \quad \quad \quad \quad \quad \quad \quad \quad \quad \quad \quad \quad \quad
  - g_i d L^{-1}  \left (  - \left (  d \omega \right )_{\vert_{ {\prescript{}{}\bigwedge}^{n+1} \mathfrak{h}_1 }}   \right ) 
  - d g_i \wedge \beta \bigg ]  \wedge \theta\\
&= 
\bigg [
-  W_{2n+1} g_i  \omega' +
d  L^{-1}  \left (   - \left (  dg_i \wedge \omega'  + g_i d \omega' \right )_{\vert_{ {\prescript{}{}\bigwedge}^{n+1} \mathfrak{h}_1 }}   \right )  \\
& \quad \quad \quad \quad \quad \quad \quad \quad \quad \quad \quad \quad \quad \quad \quad \quad \quad \quad \quad  
 - g_i d L^{-1}  \left (  - \left (  d \omega \right )_{\vert_{ {\prescript{}{}\bigwedge}^{n+1} \mathfrak{h}_1 }}   \right )  \bigg ]  \wedge \theta\\
&=
 \bigg [
- W_{2n+1} g_i  \omega' +
d  L^{-1}  \left (
   - \left (  dg_i \wedge \omega'  \right )_{\vert_{ {\prescript{}{}\bigwedge}^{n+1} \mathfrak{h}_1 }} 
   - \left (  g_i d \omega' \right )_{\vert_{ {\prescript{}{}\bigwedge}^{n+1} \mathfrak{h}_1 }} 
  \right ) \\
& \quad \quad \quad \quad \quad \quad \quad \quad \quad \quad \quad \quad \quad \quad \quad \quad \quad \quad \quad  
 - g_i d L^{-1}  \left (  - \left (  d \omega \right )_{\vert_{ {\prescript{}{}\bigwedge}^{n+1} \mathfrak{h}_1 }}   \right )  \bigg ]  \wedge \theta,\\
& \text{since $d$ and $L^{-1} $ are linear},\\
&=
\bigg [
- W_{2n+1} g_i  \omega' 
- d  L^{-1}  \left (
   \left (  dg_i \wedge \omega'  \right )_{\vert_{ {\prescript{}{}\bigwedge}^{n+1} \mathfrak{h}_1 }} 
  \right )  
- d  L^{-1}  \left (
   \left (  g_i d \omega' \right )_{\vert_{ {\prescript{}{}\bigwedge}^{n+1} \mathfrak{h}_1 }} 
  \right ) \\
& \quad \quad \quad \quad \quad \quad \quad \quad \quad \quad \quad \quad \quad \quad \quad \quad \quad \quad \quad  
+ g_i d L^{-1}  \left (  \left (  d \omega \right )_{\vert_{ {\prescript{}{}\bigwedge}^{n+1} \mathfrak{h}_1 }}   \right )  \bigg ]  \wedge \theta\\
&=
\bigg [
- W_{2n+1} g_i  \omega' 
- d  L^{-1}  \left (
    \left (  dg_i \wedge \omega'  \right )_{\vert_{ {\prescript{}{}\bigwedge}^{n+1} \mathfrak{h}_1 }} 
  \right )   
- d  \left ( g_i  L^{-1}  \left (
   \left (  d \omega' \right )_{\vert_{ {\prescript{}{}\bigwedge}^{n+1} \mathfrak{h}_1 }} 
  \right )   \right ) \\
& \quad \quad \quad \quad \quad \quad \quad \quad \quad \quad \quad \quad \quad \quad \quad \quad \quad \quad \quad  
+ g_i d L^{-1}  \left (  \left (  d \omega \right )_{\vert_{ {\prescript{}{}\bigwedge}^{n+1} \mathfrak{h}_1 }}   \right )  \bigg ]  \wedge \theta\\
&=
\bigg [
- W_{2n+1} g_i  \omega' 
- d  L^{-1}  \left (
   \left (  dg_i \wedge \omega'  \right )_{\vert_{ {\prescript{}{}\bigwedge}^{n+1} \mathfrak{h}_1 }} 
  \right )   
- d    g_i  \wedge  L^{-1}  \left (
   \left (  d \omega' \right )_{\vert_{ {\prescript{}{}\bigwedge}^{n+1} \mathfrak{h}_1 }} 
  \right )    \\
& \quad \quad \quad \quad \quad \quad \quad \quad \quad \quad \quad
-   g_i   d   L^{-1}  \left (
   \left (  d \omega' \right )_{\vert_{ {\prescript{}{}\bigwedge}^{n+1} \mathfrak{h}_1 }} 
  \right )  
+ g_i d L^{-1}  \left (  \left (  d \omega \right )_{\vert_{ {\prescript{}{}\bigwedge}^{n+1} \mathfrak{h}_1 }}   \right )  \bigg ]  \wedge \theta\\
&=
 \bigg [
- W_{2n+1} g_i   \omega' 
- d  L^{-1}  \left (
   \left (  dg_i \wedge \omega'  \right )_{\vert_{ {\prescript{}{}\bigwedge}^{n+1} \mathfrak{h}_1 }} 
  \right )  
- d    g_i  \wedge  L^{-1}  \left (
   \left (  d \omega' \right )_{\vert_{ {\prescript{}{}\bigwedge}^{n+1} \mathfrak{h}_1 }} 
  \right )     \bigg ] 
 \wedge \theta.
\end{align*}

To recap, we rewrote equation (\ref{2piecesT}) as
\begin{align}
\begin{aligned}\label{2piecesTrewritten}
&\left [ \left ( \partial G \right ) \mres{g_i} -  \partial \left ( G \mres{g_i} \right ) \right ]  ([\omega ]_{I^{n}})\\
&=
G  \Bigg (
 \bigg [
- W_{2n+1} g_i   \omega' 
- d  L^{-1}  \left (
   \left (  dg_i \wedge \omega'  \right )_{\vert_{ {\prescript{}{}\bigwedge}^{n+1} \mathfrak{h}_1 }} 
  \right )  
- d    g_i  \wedge  L^{-1}  \left (
   \left (  d \omega' \right )_{\vert_{ {\prescript{}{}\bigwedge}^{n+1} \mathfrak{h}_1 }} 
  \right )     \bigg ] 
 \wedge \theta \Bigg ).
\end{aligned}
\end{align}

This allows us to rewrite Lemma 3.11\cite{GCslicing} as
\begin{align}
\begin{aligned}\label{lastTrewritten}
&\left [ \left ( \partial G \right ) \mres{g_i} -  \partial \left ( G \mres{g_i} \right ) \right ]  (\omega)  \\
&=
\begin{cases}
G   \left ( \left [ d^{(1)} g_i \wedge  \omega \right ]_{I^{k+1}}  \right ) , 
\quad [\omega ]_{I^{k}} \in \mathcal{D}^k_{\mathbb{H}}(U) = \frac{\Omega^k}{I^k},
\quad \text{if }k<n,\\
G  \Bigg ( 
\bigg [
- W_{2n+1} g_i   \omega' 
- d  L^{-1}  \left (
   \left (  dg_i \wedge \omega'  \right )_{\vert_{ {\prescript{}{}\bigwedge}^{n+1} \mathfrak{h}_1 }} 
  \right )  
- d    g_i  \wedge  L^{-1}  \left (
   \left (  d \omega' \right )_{\vert_{ {\prescript{}{}\bigwedge}^{n+1} \mathfrak{h}_1 }} 
  \right )     \bigg ] 
 \wedge \theta
\Bigg ) , \\
\quad \quad  \quad  \quad \quad  \quad  \quad  \quad \quad  \quad  \quad \quad  \quad  \quad \quad  \quad  \quad [\omega ]_{I^{n}} \in \mathcal{D}^n_{\mathbb{H}}(U) = \frac{\Omega^n}{I^n},
\quad \text{if }k=n,\\
G   \left ( \left ( d^{(1)} g_i \wedge  \omega \right )_{\vert_{J^{k+1} }}  \right ) , 
\quad \omega \in \mathcal{D}^k_{\mathbb{H}}(U) = J^k,
\quad \text{if }k>n.
\end{cases}
\end{aligned}
\end{align}


\section{Constraints and the special case of $\mathbb{H}^1$}

\noindent
We would now like to remind the reader that our goal was to find an inequality similar to property (4) in Proposition \ref{next3properties} for the case $k=n$. Specifically, we tried to find a result equivalent to Lemma 3.12\cite{GCslicing}, for the case $k=n$. We were unable to achieve this and now we will show some of the difficulties in following the idea of the proof of 3.12\cite{GCslicing}, for the special case of $\mathbb{H}^1$.\\

\noindent
As an aside, it is important to notice that the case $n=1$ for some aspects behaves differently from $n\geq 2$ (for more details, see the note after Proposition 3.6 in \cite{GCslicing}). Hence these difficulties should not be generalised, a priori, to all of $\mathbb{H}^n$.\\

\noindent
In $\mathbb{H}^1$, we consider $ [\omega ]_{I^{1}} \in \mathcal{D}^1_{\mathbb{H}}(U) $ and we can take $\omega = \omega' =\omega_1 dx + \omega_2 dy$. In the case $k=n$ of equation (\ref{lastTrewritten}), the input to the current $G$ is a sum of three components. The first one is
\begin{align*}
- T g_i \omega' \wedge \theta =& T g_i \omega_1 dx \wedge \theta + T g_i \omega_2 dy \wedge \theta.
\end{align*}
Now we look at the second component. Since
$$
d g_i = X g_i dx + Y g_i dy + T g_i \theta \quad \text{and} \quad \omega' = \omega_1 dx + \omega_2 dy,
$$
we have
$$
d g_i \wedge \omega' = X g_i \omega_2 dx \wedge dy + Y g_i \omega_1 dy \wedge dx + T g_i \omega_1 \theta \wedge dx + T g_i \omega_2 \theta \wedge dy.
$$
Then
\begin{align*}
\left ( d g_i \wedge \omega' \right )_{\vert_{ {\prescript{}{}\bigwedge}^{2} \mathfrak{h}_1 }}  =& \left ( X g_i \omega_2 - Y g_i \omega_1 \right ) dx \wedge dy ;\\
\left ( d g_i \wedge \omega' \right )_{\vert_{ {\prescript{}{}\bigwedge}^{2} \mathfrak{h}_1 }}  =& - \left ( X g_i \omega_2 - Y g_i \omega_1 \right ) d \theta ;\\
- L^{-1}   \left ( \left ( d g_i \wedge \omega' \right )_{\vert_{ {\prescript{}{}\bigwedge}^{2} \mathfrak{h}_1 }}  \right ) =&  X g_i \omega_2 - Y g_i \omega_1;\\
- dL^{-1}   \left ( \left ( d g_i \wedge \omega' \right )_{\vert_{ {\prescript{}{}\bigwedge}^{2} \mathfrak{h}_1 }} \right ) =& X \left (  X g_i \omega_2 - Y g_i \omega_1 \right ) dx + Y \left (  X g_i \omega_2 - Y g_i \omega_1 \right ) dy + T \left (  X g_i \omega_2 - Y g_i \omega_1 \right ) \theta;\\
- dL^{-1}   \left ( \left ( d g_i \wedge \omega' \right )_{\vert_{ {\prescript{}{}\bigwedge}^{2} \mathfrak{h}_1 }} \right ) \wedge \theta =& X \left (  X g_i \omega_2 - Y g_i \omega_1 \right ) dx \wedge \theta + Y \left (  X g_i \omega_2 - Y g_i \omega_1 \right ) dy \wedge \theta.
\end{align*}
Turning to the third component, we have:
\begin{align*}
d \omega' =& Y\omega_1 dy \wedge dx + T\omega_1 \theta \wedge dx     + X\omega_2 dx \wedge dy + T\omega_2 \theta \wedge dy;\\
\left ( d \omega' \right )_{\vert_{ {\prescript{}{}\bigwedge}^{2} \mathfrak{h}_1 }}  =& \left (  X\omega_2 - Y\omega_1 \right ) dx \wedge dy ;\\
\left ( d \omega' \right )_{\vert_{ {\prescript{}{}\bigwedge}^{2} \mathfrak{h}_1 }}  =& - \left (  X\omega_2 - Y\omega_1 \right ) d \theta ;\\
- L^{-1}   \left ( \left ( d \omega \right )_{\vert_{ {\prescript{}{}\bigwedge}^{2} \mathfrak{h}_1 }}  \right ) =&  X\omega_2 - Y\omega_1. 
\end{align*}
Thus, since $d g_i = X g_i dx + Y g_i dy + T g_i \theta$, we have
$$
- dg_i \wedge L^{-1}   \left ( - \left ( d \omega \right )_{\vert_{ {\prescript{}{}\bigwedge}^{2} \mathfrak{h}_1 }}  \right )    \wedge  \theta = \left (  X\omega_2 - Y\omega_1 \right )  X g_i dx \wedge \theta +
\left (  X\omega_2 - Y\omega_1 \right )  Y g_i dy \wedge \theta.
$$
Combining the three components, we get
\begin{align*}
&\bigg [
- W_{2n+1} g_i   \omega' 
- d  L^{-1}  \left (
   \left (  dg_i \wedge \omega'  \right )_{\vert_{ {\prescript{}{}\bigwedge}^{n+1} \mathfrak{h}_1 }} 
  \right )  
- d    g_i  \wedge  L^{-1}  \left (
   \left (  d \omega' \right )_{\vert_{ {\prescript{}{}\bigwedge}^{n+1} \mathfrak{h}_1 }} 
  \right )     \bigg ] 
 \wedge \theta\\
&=T g_i \omega_1 dx \wedge \theta + T g_i \omega_2 dy \wedge \theta \\
&+ X \left (  X g_i \omega_2 - Y g_i \omega_1 \right ) dx \wedge \theta + Y \left (  X g_i \omega_2 - Y g_i \omega_1 \right ) dy \wedge \theta \\
&+ \left (  X\omega_2 - Y\omega_1 \right ) X g_i dx \wedge \theta + 
\left (  X\omega_2 - Y\omega_1 \right )  Y g_i dy \wedge \theta \\
&= \big ( T g_i \omega_1 + X \left (  X g_i \omega_2 - Y g_i \omega_1 \right ) + \left (  X\omega_2 - Y\omega_1 \right ) X g_i \big ) dx \wedge \theta\\
&+ \big ( T g_i \omega_2 + Y \left (  X g_i \omega_2 - Y g_i \omega_1 \right ) + \left (  X\omega_2 - Y\omega_1 \right )  Y g_i \big ) dy \wedge \theta.
\end{align*}
Ideally we would like to bound these terms with constants, but this does not seem possible. %
The tools used in article \cite{GCslicing}, %
are insufficient in this situation, as they cannot bound the first derivatives of the coefficients of $\omega$ as well as the second derivatives of the approximating functions $g_i$.\\

\noindent
In conclusion, it seems that to obtain results in the case $k=n$ with this approach, one would need to restrict some hypotheses. Other approaches that we did not explore deeply are
\begin{enumerate}
\item
consider the differences between $\mathbb{H}^1$, $\mathbb{H}^2$ and $n \geq 2$ in general (see the note after Proposition 3.6 in article \cite{GCslicing}),
\item
use the local definition of currents (see Definition 2.3 in article \cite{LANG}),
\item
use currents with locally finite mass (see end of section 2 in  article \cite{GCslicing} and, mainly, Proposition 5.15 and Definition 5.18 in article \cite{FSSC}),
\item
deal with the second derivative of $g_i$ by restricting how it approximates the function $f$ (see Lemma 3.9 in article \cite{GCslicing} and 4.3.1 in article \cite{FED}).
\end{enumerate}



\bibliography{Bibliography_250822} 
\bibliographystyle{abbrv}

\end{document}